\documentclass{article}

\usepackage{arxiv}

\usepackage[utf8]{inputenc} 
\usepackage[T1]{fontenc}    
\usepackage{hyperref}       
\usepackage{url}            
\usepackage{booktabs}       
\usepackage{amsfonts}       
\usepackage{nicefrac}       
\usepackage{microtype}      
\usepackage{lipsum}
\usepackage{amsmath}
\usepackage{graphicx}
\usepackage{enumerate}
\usepackage{lineno}
\usepackage{setspace}
\title{Optimising collective accuracy among rational individuals in sequential decision-making with competition}

\author{
  Richard P. Mann \\
  Department of Statistics, School of Mathematics, University of Leeds\\
  \texttt{r.p.mann@leeds.ac.uk}
}

\begin{document}
\maketitle

\begin{abstract}
Theoretical results underpinning the Wisdom of Crowds, such as the Condorcet Jury Theorem, point to substantial accuracy gains through aggregation of decisions or opinions, but the foundations of this theorem are routinely undermined in circumstances where individuals are able to adapt their own choices based after observing what other agents have chosen. In sequential decision-making, rational agents use the choices of others as a source of information about the correct decision, creating powerful correlations between different agents' choices that violate the assumptions of independence on which the Condorcet Jury Theorem depends. In this paper I show how such correlations emerge when agents are rewarded solely based on their individual accuracy, and the impact of this on collective accuracy. I then demonstrate how a simple competitive reward scheme, where agents' rewards are greater if they correctly choose options that few have already chosen, can induce rational agents to make independent choices, returning the group to optimal levels of collective accuracy. I further show that this reward scheme is robust, offering improvements to collective accuracy across of wide range of competition strengths, suggesting that such schemes could be effectively implemented in real-world contexts to improve collective wisdom.
\end{abstract}

\section*{Introduction}
It has long been recognised that aggregating the opinions, estimates or decisions of many individuals can give superior results compared to relying on a single individual alone \cite{de2014essai, galton1907vox}. Sometimes termed The Wisdom of Crowds \cite{surowiecki2005wisdom}, such aggregation is a simple but potentially powerful example of collective intelligence, and one that acts as both a justification for democratic decision-making institutions \cite{list2001epistemic, landemore2012democratic} and a motivation for utilising fora such as social media to harness the potential of global collective knowledge \cite{klein2011harvest}.

The Condorcet Jury Theorem (CJT) \cite{de2014essai, boland1989majority} demonstrates that collective accuracy, in the form of a majority vote, can far exceed individual accuracy under an idealised assumption that agents choose independently. While the CJT has motivated many appeals to The Wisdom of Crowds (e.g. \cite{list2001epistemic,list2004democracy, king2007use}), in reality this independence assumption is routinely violated in collective decision-making scenarios where agents are able to observe each other and use social information to motivate their own choices \cite{surowiecki2005wisdom}. The Wisdom of Crowds requires that a group must effectively aggregate the private information held by its members, but information cascades can result from social learning, such that within a group a large proportion of individuals simply follow the decisions made by others, without reference to any private information they may have \cite{bikhchandani1992theory}. Empirical studies have demonstrated how readily humans copy the actions of others \cite{faria2010coll,gallup2012visual, mann2013tdo}, in common with other animals \cite{sumpter2009quo}, when those actions are readily observable. The tendency of agents to follow the decisions of others can be rational from an individual perspective \cite{bikhchandani1992theory, anderson1997information,mann2018cdm, tump2020wise}, but such self-reinforcing cascades of social information can cause very large scale errors in collective judgement, as illustrated anecdotally in the historical examples given by Mackay in `Extraordinary Popular Delusions and the Madness of Crowds' \cite{mackay2012extraordinary}. Scientific study also suggests that, under controlled conditions, allowing individuals to update their own beliefs in the light of observing others tends to reduce the accuracy of collective estimations \cite{lorenz2011hsi}, even as it increases the average accuracy of individual agents \cite{tump2020wise}.

The dangers of relying on social information are highly pertinent since sequential decision-making is common across a wide variety of domains. We often choose what to buy, where to eat, or even how to vote based on the choices or expressed opinions of others before us. Sequential decisions may be present even when a system is designed to elicit individuals' independent decisions. Consider for example the case of formalised peer-review of scientific publications or grant proposals. Here, reviewers apparently provide their reviews independently, but this ignores the effect of author status, which provides a proxy for the decisions of past reviewers of the same author. Likewise, characteristics such as the fame of an individual, or the market share of a product, may serve to indicate a preponderance of past choices made, even if these are not directly observed. Advertising that points to the number of users or consumers of a product is suggestive of the influence such past decisions can have on future purchases.

Given the prevalence of sequential decision-making across many areas in which we may wish to access collective knowledge, how might we overcome its deleterious effects upon collective wisdom? One potential solution is the introduction of competition between agents who make the same choice \cite{hong, mann2017optimal}, thus penalising agents who follow others. Previous work on sequential decision-making \cite{perez2011cab, arganda2013acr, mann2018cdm, mann2020collective} has assumed that rewards are independent of which choices other agents make, with such choices being useful only as a source of information about the rewards available in the environment. In this paper I extend this framework to allow for rewards that depend intrinsically upon the choices made by others, such that an option may become more or less rewarding based on how many other agents have also chosen it. Using this model, I show how social information in the absence of competition can reduce the collective accuracy of a group, and how introducing competition in the form of diminishing rewards for options already chosen by other agents can eliminate correlations between agents' choices, and return the group to an optimal level of collective accuracy.

\section*{Model}
I consider a binary choice scenario with potential options labelled as A and B. In any given choice, one option is `correct' and the other is `incorrect'. This scenario is similar to that in \cite{mann2018cdm}, and the model described below largely follows the framework developed in that paper. 

Agents sequentially choose either A or B, and are able to observe all choices made before their own, such that these such constitute common knowledge \cite{aumann1976agreeing}. The choice made by individual $i$ can be labelled as $C_i= 1$ if A is chosen, or $C_i=-1$ if B is chosen, and a sequence of $k$ decisions $S$ is an ordered series $C_1, C_2, \ldots C_k$. The collective decision is defined as the majority choice when all $n$ agents have decided, and for simplicity I consider only cases where $n$ is odd so there are no tied collective decisions.

Agents choose between the two options based on reward criteria and their own inferences about the probability that each option is correct, so as to maximise their expected reward. The true state of the world can be given by a variable $x$, which takes the value $x=1$ if option A is correct, and $x=-1$ if option B is correct. All agents are assumed to share a common, symmetric and uninformative prior about the value of $x$:
\begin{equation}
    P(x=1) = P(x=-1) = 1/2.
\end{equation}
Agents are informed by two sources of information. The first is a noisy private signal $\Delta_i$ received independently by each agent $i$, with variance $\epsilon$:
\begin{equation}
    p(\Delta_i \mid x) = \frac{1}{\epsilon}\phi\left(\frac{\Delta_i-x}{\epsilon}\right),
\end{equation}
where $\phi(\cdot)$ is the standardised normal probability distribution function.
The second source of information is the social information provided by the sequence of previous decisions $S$. Agents update their knowledge of $x$ by performing Bayesian inference:
\begin{equation}
    P(x \mid \Delta_i, S) = \frac{P(x)P(S \mid x)p(\Delta_i \mid x)}{\sum_{x \in \{-1,1\}} P(x)P(S \mid x)p(\Delta_i \mid x)}
\end{equation}
where the equation above makes use of the assumed independence of private signals, and thus the independence of $S$ and $\Delta_i$ conditioned on $x$.

\subsection*{Rewards}
Agents are motivated to make accurate choices by a retrospective reward policy that assigns rewards once the true correct choice is known. A simple and intuitive policy is to reward agents if they made the correct choice, thus motivating each individual to be as accurate as possible. This can be labelled as `binary' rewards, in common with previous models on simultaneous decision-making \cite{hong, mann2017optimal}, since agents receive a reward of either zero or one (in some standardised reward units) for each choice. This reward policy can be defined mathematically via a reward function $r(C_i, S, x)$ that depends on the choice, $C_i$, made by individual $i$, the sequence of past decisions $S$ and the true state of the world $x$, with binary rewards being defined as:
\begin{equation}
    r_{\textrm{binary}}(C_i, S, x) = \delta_{C_i, x},
\end{equation}
where $\delta_{l,k}$ is the Kronecker delta function.

A binary reward function attributes rewards based solely on the accuracy of an individual's choice, and is independent of the decisions made by others. More generally we can consider a reward scheme that depends on past decisions that the agent can observe. A simple way to do this is to make modulate the reward with a function that depends on the individual choice and $S$:
\begin{equation}
    r_{\textrm{general}}(C_i, S, x) = f(C_i, S)\delta_{C_i, x}.
\end{equation}
This continues to reward (and thus incentivise) accuracy through the $\delta_{C_i, x}$ term, but can also directly reward or penalise choosing the same option as others, thus incentivising either conformity or diversity of choices.
\subsection*{Rational individual choice}
Given a reward function $r(C_i, S, x)=f(C_i, S)\delta_{C_i, x}$, an agent can evaluate the expected reward $\mathbb{E}(r_A \mid S, \Delta_i)$ from choosing A, conditioned on the available private and social information:
\begin{equation}
\begin{split}
    \mathbb{E}(r_A \mid S, \Delta_i) &= r(C_i=1, S, x=1)P(x=1 \mid \Delta_i, S) + r(C_i=1, S, x=-1)P(x=-1 \mid \Delta_i, S) \\
    &= f(C_i=1, S)P(x=1 \mid \Delta_i, S)
    \end{split}
\end{equation}
and similarly for choosing B:
\begin{equation}
    \mathbb{E}(r_B \mid S, \Delta_i) = f(C_i=-1, S)P(x=-1 \mid \Delta_i, S), 
\end{equation}
According to the principle of expected reward maximisation, a rational agent will then select A if and only if $\mathbb{E}(r_A \mid S, \Delta_i > \mathbb{E}(r_B \mid S, \Delta_i)$ (since $\Delta_i$ is real-valued, a tied expectation has zero probability mass). Using the general reward function above, this condition simplifies to:
\begin{equation}
    \mathbb{E}(r_A \mid S, \Delta_i) > \mathbb{E}(r_B \mid S, \Delta_i) \Rightarrow \frac{P(x=1 \mid \Delta_i, S)}{P(x=-1 \mid \Delta_i, S)} > \frac{f(C_i = -1, S)}{f(C_i = 1, S)}.
\end{equation}
That is, for an agent to choose A, its assessment of the difference in probability for A to be correct rather than B must outweigh any penalty it receives for choosing A over B based on the past decisions.

A feature of the above decision-making procedure is that there exists some critical value of an agent's private information, $\Delta_i^*$, which would make the expected reward of choosing A or B equal:
\begin{equation}
    \mathbb{E}(r_A \mid S, \Delta_i^*) = \mathbb{E}(r_B \mid S, \Delta_i^*).
\end{equation}
This implies that agent $i$ will choose A if and only if $\Delta_i > \Delta_i^*$. Substituting the definition of the expected reward and the conditional probability $P(x \mid \Delta_i, S)$, this gives:
\begin{equation}
    P(S \mid x=1)\phi((\Delta_i^*-1)/\epsilon)f(C_i = 1, S) = P(S \mid x=-1)\phi((\Delta_i^*+1)/\epsilon)f(C_i = -1, S).
\end{equation}
We can recognise that $\frac{\phi((\Delta_i^*-1)/\epsilon)}{\phi((\Delta_i^*+1)/\epsilon)} = \exp(2\Delta_i^*/\epsilon^2)$, and thus the expression above can be rearranged to give:
\begin{equation}
    \Delta_i^* = \frac{\epsilon^2}{2}\left( \log\left( \frac{P(S \mid x=-1)}{P(S \mid x=1} \right) + \log\left( \frac{f(C_i = -1, S)}{f(C_i = 1, S)} \right) \right)
    \label{eqn:threshold}
\end{equation}
Since subsequent agents are able to observe the value of $S$ that agent $i$ was responding to, they can calculate the corresponding value of $\Delta_i^*$. Combined with observing the decision agent $i$ makes, this enables them to infer whether agent $i$'s private information was greater than or less than this threshold value. The probability of a sequence $S$, conditioned on $x$, can therefore be evaluated with reference to each of the thresholds calculated for the previous agents:
\begin{equation}
\begin{split}
    P(S \mid x) &= \prod_{j \in \textrm {choosing A}} P(\Delta_j > \Delta_j^*) \times \prod_{j \in \textrm {choosing B}} P(\Delta_j < \Delta_j^*) \\
    &= \prod_{j < i} \Phi(C_j(x-\Delta_j^*)/\epsilon)
    \end{split}
    \label{eqn:Psx}
\end{equation}
where $\Phi(\cdot)$ is the cumulative probability function of the standard normal distribution. Since the thresholds depend themselves on the past sequence of decisions, the probability of a sequence can be evaluated recursively by calculating the threshold for each sub-sequence.

\section*{Results}
\subsection*{Social response under binary rewards}
\begin{figure}
    \centering
    \includegraphics[width=15cm]{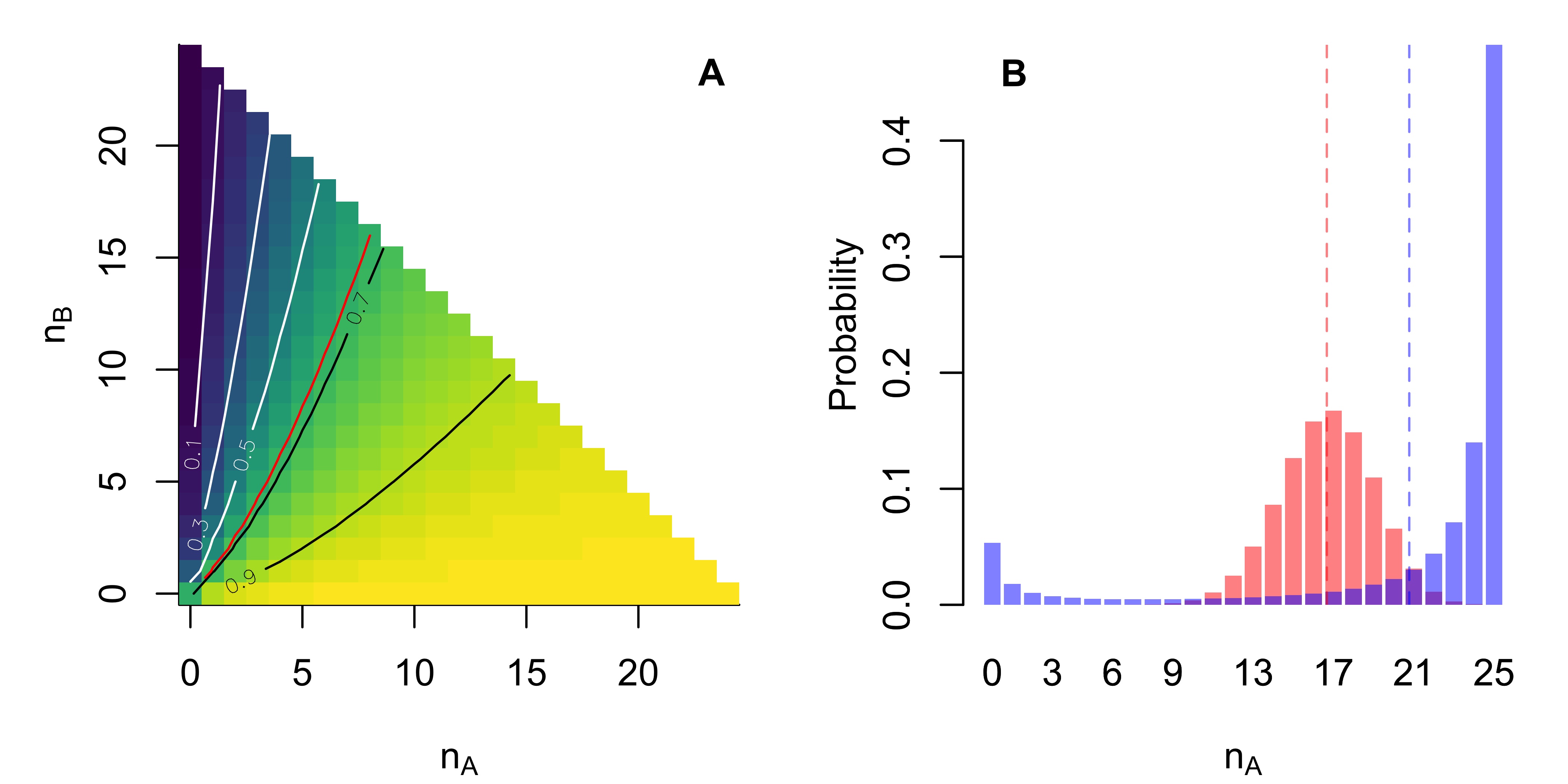}
    \caption{Characterising the response to social information in sequential decisions under binary rewards. (A) the probability that an agent will choose option A when that is the correct choice, conditioned on the number of previous decisions for options A and B, averaging over all sequences consistent with those aggregate number of decisions. In this example the environmental noise is $\epsilon = 2.32$, giving an individual choice accuracy of $q=2/3$; the red contour line indicates this probability. (B) The probability for $n_A$ agents to select option A when that is the correct choice, averaged over a full sequence of decisions in a group of $n=25$ agents. The blue bars indicate the probability when rational agents are subject to binary rewards, red bars indicate the probability if all agents select independently. The dashed lines indicate the mean of each probability distribution. Agents responding rationally to binary rewards have a higher average number of individually-successful decisions, but a lower probability of a correct majority decision.}
    \label{fig:socialresponse}
\end{figure}

The influence of social information on decision making can be characterised by observing its effect on both individual decisions and on the aggregate outcomes in groups. A simple way to visualise the influence of previous choices on an individual's decision is via the probability that a focal individual will choose the correct option, arbitrarily taken to be option A, conditioned on there having previously been $n_A$ and $n_B$ agents choosing A and B respectively. This probability is shown in Figure \ref{fig:socialresponse}A, assuming that agents are responding to binary rewards ($f(C_i = 1, S) = f(C_i = -1, S) = 1$). Because the decision to choose A or B depends in theory on both the full sequence of previous choices and the agent's private information, the probability shown in this figure is a weighted average over all sequences consistent with specified values of $n_A$ and $n_B$, and all possible values of the focal agent's private information:
\begin{equation}
    P(C_i=1 \mid n_A, n_B) = \sum_{S \in n_A,n_B} P(\Delta_i > \Delta_S^*)P(S \mid x=1),
\end{equation}
where the summation is over the set of all sequences with $n_A$ and $n_B$ individuals choosing A and B. 

This figure shows that agents respond strongly to the decisions made by others, such that the probability to choose A is highly dependent on the values of $n_A$ and $n_B$. In particular, in most cases where $n_B > n_A$ the focal agent is less likely to choose the correct option than they would be if they chose independently; the red contour lines indicates this independent choice probability. This implies that incorrect decisions by agents at the beginning of the sequence can lead to a cascade of later agents also making incorrect choices. This is reflected in distribution of aggregate outcomes at the group level, characterised by the probability that $n_A$ agents will select option A in total. This is plotted in Figure \ref{fig:socialresponse}B for both the case of independent decisions (red bars) and for agents using social information with binary rewards (blue bars). This plot shows the dramatic difference in aggregate outcomes that results from social information use. When agents choose independently the aggregate outcomes are clustered in a binomial distribution that peaks at the the mean value of $n_A= n\Phi(1/\epsilon)$, with a very low probability that fewer than half the agents choose A. Under social information the aggregate outcomes become bimodal, with a large peak at $n_A=n$ and a secondary peak at $n_A=0$. The result of this is that the mean number of correct decisions increases (compare the blue and red dashed lines), but there is a much greater probability that a majority of agents will choose the incorrect option (B). As such, each individual is more likely to choose the correct option, but the majority choice of the group is less likely to be correct.

\subsection*{Condorcet-retrieving reward function}
Under binary rewards, agents tend to follow past decisions with increasing strength over the course of a sequence of choices (Figure \ref{fig:socialresponse}A). Since this breaks the assumption of independence in the CJT, it also reduces the accuracy of collective decisions as defined by the majority choice, as shown in Figure \ref{fig:socialresponse}B. To improve collective accuracy it is therefore necessary to reduce the correlation between decisions. If agents make choices independently, this implies that the threshold value of $\Delta_S^*$ should be independent of the value of $S$. Since agents begin with a symmetric prior $P(x=1) = 1/2$, it further implies that this threshold must be zero - i.e. agents will choose A or B based solely on the direction of their private information. One can therefore retrieve independent choices, and thus the accuracy implied by the CJT, by seeking a reward function $r_{\textrm{condorcet}}(C_i, S, x)$ such that:
\begin{equation}
    \mathbb{E}(r_A \mid S, \Delta_i=0) = \mathbb{E}(r_B \mid S, \Delta_i=0) \ \forall \ S.
\end{equation}
Expanding the definition for the expected reward, this gives:
\begin{equation}
        P(x=1 \mid \Delta_i=0, S)f(C_i=1, S) = P(x=-1 \mid \Delta_i=0, S)f(C_i = -1, S).
\end{equation}
By substituting Bayes rule for the conditional probability of $x$, we get:
\begin{equation}
        \frac{p(\Delta_i=0 \mid x=1)P(S \mid x=1)}{ p(\Delta_i=0 \mid x=-1)P(S \mid x=-1)} = \frac{f(C_i=-1, S)}{f(C_i=1, S}
\end{equation}
By construction, under this Condorcet reward scheme, all thresholds for private information are zero. The probability $P(S \mid x)$ thus simplifies to a product of independent choices:
\begin{equation}
    P(S \mid x) = \Phi(x/\epsilon)^a\Phi(-x/\epsilon)^b,
\end{equation}
where $a$ is the number of agents who have previously chosen A and $b$ the number who have chosen B. Substituting this expression and recognising that $p(\Delta_i=0 \mid x=1) = p(\Delta_i=0 \mid x=-1)$ , we therefore get:
\begin{equation}
        \frac{\Phi(1/\epsilon)^a\Phi(-1/\epsilon)^b}{\Phi(-1/\epsilon)^a\Phi(1/\epsilon)^b} = \frac{f(C_i=-1, S)}{f(C_i=1, S}
\end{equation}
This expression can be simplified by defining $q = \Phi(1/\epsilon)$ as the probability that a single agent will independently choose the correct option. This then reduces to:
\begin{equation}
\begin{split}
     \frac{f(C_i=-1, S)}{f(C_i=1, S} &= \frac{q^a(1-q)^b}{(1-q)^a q^b} \\
     &= Q^{a-b},
    \end{split}
\end{equation}
where $Q = q/(1-q)$. This expression can be satisfied by a reward scheme:
\begin{equation}
    r_{\textrm{condorcet}}(C_i, S, x)= \begin{cases}
    \delta_{C_i, x}Q^{-a} \ \textrm{if} \ C_i = 1\\
    \delta_{C_i, x}Q^{-b} \ \textrm{if} \ C_i = -1
    \end{cases}
\end{equation}
This expression shows that rational agents can be motivated to make independent choices if the rewards for each option are reduced geometrically with the number of agents that have already chosen that option. This is a very convenient reward system for several reasons: First, it is symmetric in the way it treats both options, so neither option needs to be arbitrarily favoured or penalised. Second, the form of the required penalty for each option depends only on the number of agents that have previous chosen it, so these penalties can be implemented locally without reference to the number choosing the other option, or the order in which those choice were made. Third, it resembles a form of competition, with each agent exhausting a fixed proportion of the potential reward remaining for the option it chooses. The geometric reduction in rewards means that for any group size the total rewards available from each option are bounded by:
\begin{equation}
    \textrm{Maximum reward} = 1 + 1/Q + 1/Q^2 \ldots = Q/(Q-1).
\end{equation}
Similarly, the expected total reward can be calculated as:
\begin{equation}
\begin{split}
    \mathbb{E}(\textrm{total reward}) &= \mathbb{E}(1 + 1/Q + 1/Q^2 \ldots 1/Q^k), \ k \sim \textrm{Bin}(n, q) \\
    &= \frac{Q-2^n(1-q)^n}{Q-1}
    \end{split}
\end{equation}
Any system that assigns rewards under this scheme can therefore estimate and bound the total rewards it would potentially need to allocate. It is notable that high values of $Q$ indicate problems that are relatively simple for individual decision makers, and these represent the lowest expectation and bound on total rewards; this naturally allows a reward system to allocate the greatest reward budget to the most difficult problems. 

\subsection*{Robustness of collective accuracy under varying competition}
The reward scheme derived above is constructed so as to maximise the accuracy of the majority choice by making individual rational decisions statistically independent, and it accomplishes this through imposing a specific form of competitive penalty. As discussed above, this form of competitive penalty has many agreeable features for implementation in real world decisions problems. However, selecting the precise strength of the competitive penalty requires knowing in advance how difficult the decision problem is, i.e. knowing the value of $Q$. In general it is unlikely that this would be precisely known in advance, although a system designer may have some intuition about whether a given decision is easy or difficult. As such, it is important to assess how robust such reward system is to misspecification of the competition strength. To do this, we can evaluate the collective accuracy under a reward scheme with variable competition strength $\beta$:
\begin{equation}
    r_{\textrm{competitive}}(C_i, S, x)= \begin{cases}
    \delta_{C_i, x}\beta^{-a} \ \textrm{if} \ C_i = 1\\
    \delta_{C_i, x}\beta^{-b} \ \textrm{if} \ C_i = -1,
    \end{cases}
\end{equation}
where we know from the above argument that the optimum value of $\beta$ should be $Q$. Under this reward scheme, the relation for critical thresholds given by equation \ref{eqn:threshold} can be simplified and evaluated efficiently as:
\begin{equation}
    \Delta_i^* = \frac{\epsilon^2}{2}\left( \sum_{j < i} \left[ \log \Phi(C_j(1-\Delta_j^*)/\epsilon) - \log \Phi(C_j(-1-\Delta_j^*)/\epsilon) \right]  + (n_A-n_B)\log\beta  \right),
\end{equation}
where $n_A$ and $n_B$ are the number of decisions for A and B respectively within the sequence $S$. Recognising the recursive pattern this further simplifies to:
\begin{equation}
    \Delta_i^* = \Delta_{i-1}^*  + \frac{\epsilon^2}{2}\left (\log \Phi(C_{i-1}(1-\Delta_j^*)/\epsilon) - \log \Phi(C_{i-1}(-1-\Delta_{i-1}^*)/\epsilon)  + C_{i-1} \log\beta  \right),
\end{equation}

The expected collective accuracy under this reward scheme can be evaluated by directly calculating the expected proportion of accurate majority decisions as a function of the adjustable competition parameter $\beta$. This is done by calculating the probability of every possible sequence of decisions in a group of $n$ agents (hence $2^n$ possible sequences) for $x=1$ and summing the probability of that set of sequences where the majority of decisions are for the correct option A:
\begin{equation}
    \mathbb{E}(\textrm{collective accuracy}) = \sum_{S \in n_A>n_B} P(S \mid x=1),
\end{equation}
where $P(S \mid x=1)$ is given by evaluating equation \ref{eqn:Psx}.

\begin{figure}
    \centering
    \includegraphics[width=15cm]{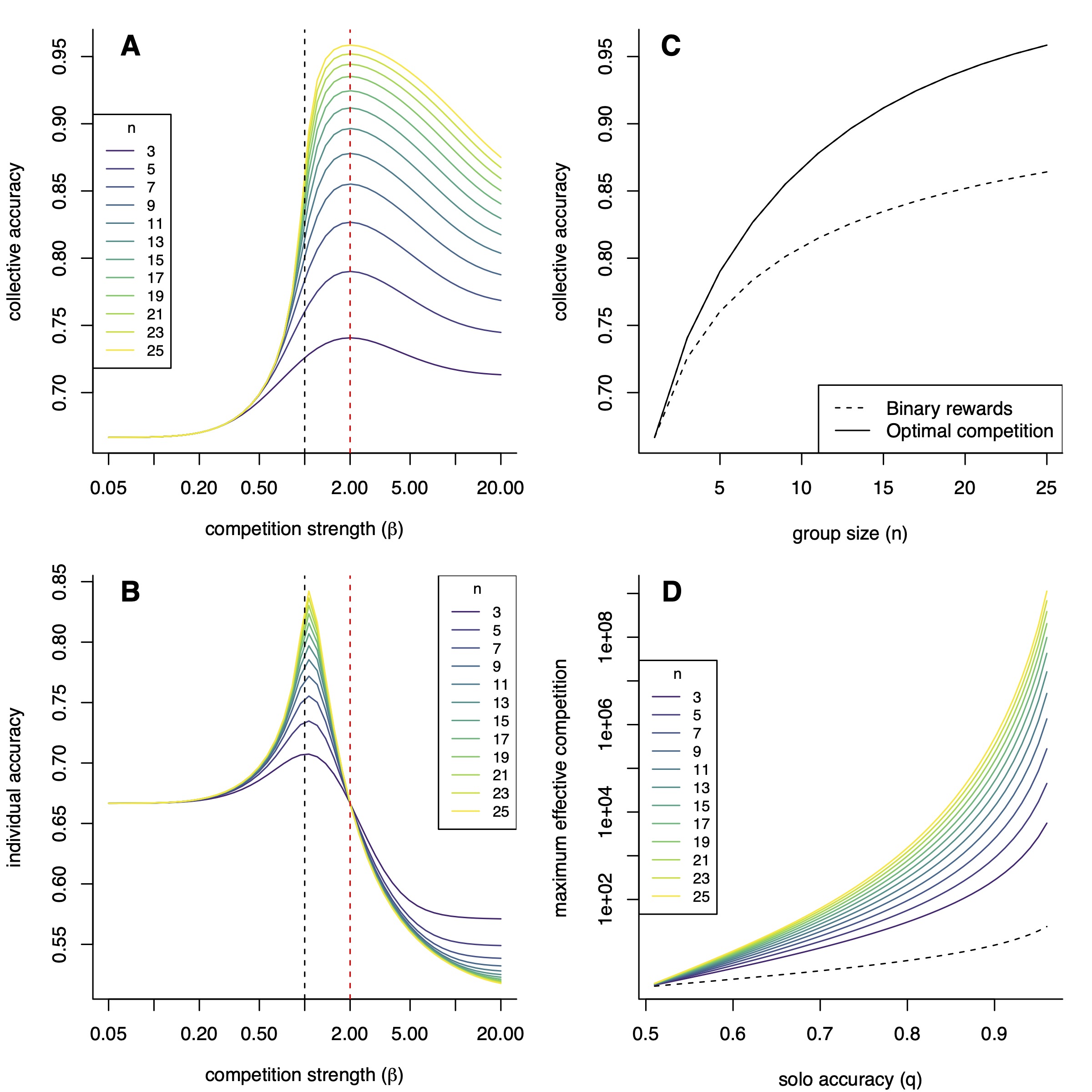}
    \caption{The effect of competition on collective accuracy. (A) With $q=2/3$, across different group sizes ($n$) collective accuracy increases with increasing competition ($\beta$) up to an optimal value given by $\beta = Q$ (indicated by the dashed red line), where collective accuracy matches that predicted by the Condorcet Jury Theorem. Higher levels of competition reduce collective accuracy, with sufficiently high values of $\beta$ leading to lower collective accuracy than under binary rewards ($\beta = 1$, indicated by the dashed black line) . Negative values of $\beta$, indicating rewards for conformity, always lead to lower collective accuracy; (B) Individual accuracy is maximised at values of $\beta$ close to one, indicating weak positive competition, and increases with group size. At the optimal competition for collective accuracy, individual accuracy is the same for all group sizes as agents choose independently; (C) The collective accuracy under optimal competition (solid line) compared to that achieved under binary rewards (dashed line) as a function of group size; (D) The maximum value of $\beta$ for which competitive rewards outperform binary rewards, for varying group size and as a function of $q$ (representing the probability for a solo agent to choose correctly). The dashed line shows the optimal value of $\beta=Q$ for comparison. The range of effective competition values (those that improve on binary rewards) is greater for easier decisions and in larger group sizes. Note the logarithmic scale on the y-axis.}
    \label{fig:competition}
\end{figure}

Figure \ref{fig:competition}A shows the collective accuracy as a function of $\beta$ for group sizes from $n=3$ to $n=25$ with an environmental noise level set of $\epsilon =2.32$, implying $q=2/3$ and $Q=2$ (i.e. individuals will make the correct choice twice as often as the wrong choice when choosing alone). This demonstrates a clear peak in accuracy in each case at the expected value of $\beta =2$, indicated by the dashed red line. At this optimum point collective accuracy matches that expected from the CJT. While a range of values of $\beta > 1$ induce greater collective accuracy than under binary rewards ($\beta =1)$, value of $\beta < 1$, which reward agents for copying past decisions, dramatically reduce collective accuracy. Figure \ref{fig:competition}B shows the individual accuracy for the same range of group sizes and competition strengths, demonstrating that increases in collective accuracy induced by competition lead to decreases in individual accuracy -- collective accuracy is maximised when individual accuracy falls to that expected from a single agent without social information, as this is when agents choose independently. Average individual accuracy is maximise at values of $\beta$ slightly greater than one, since under binary rewards each agent is motivated to maximise its own accuracy, without regard for the value of the social information it provides to those further along the sequence of decision makers (cf. \cite{torney2015social}).
Figure \ref{fig:competition}C shows the relationship between the collective accuracy achieved by the Condorcet reward scheme and that achieved without competition (binary rewards), showing that the effect is stronger in larger groups, which suffer relatively more from information cascades under binary rewards. Although collective accuracy is maximised when competition is optimised to produce independent decision making, there is a range of values of $\beta$ which induce greater collective accuracy than under binary rewards, as seen in Figure \ref{fig:competition}A. The size of this range shows how well-tuned competition must be to generate improvements in collective accuracy, and thus is indicative of how plausible effectively implementing such a reward scheme might be in practice. Figure \ref{fig:competition}D shows the maximum value of $\beta$ that outperforms binary rewards as a function of $q$, for group sizes from $n=3$ to $n=25$ (solid lines), as well as the optimal value of $\beta$ for comparison (dashed line). Inherently easier decisions permit a greater range of effective competition strengths, and this range increases very rapidly as $q$ approaches one (note the log scale on the y-axis). Larger groups also permit a wider range of effective competition strengths, even though the optimal competition strength does not depend on $n$.

\section*{Discussion}
When agents are rewarded solely for their individual accuracy they tend to follow previous decisions. While this increases the expected proportion of agents that make the correct choice, it reduces the probability that the majority of agents is correct compared to agents who make their decisions independently. Errors in early decisions can make subsequent decision-makers less accurate than they would have been alone. Hence, while on average individually beneficial, social information is deleterious to anyone seeking to use the Wisdom of Crowds by relying on the majority opinion.  

Social information may potentially be restricted exogenously, by insisting that individuals make their choices without access to the choices made by others. However, such a scenario requires tight control of the information individuals have access to, and is unlikely to be plausible when making use of collective wisdom in real-world contexts such as online review systems and social media \cite{klein2011harvest}. 

Here I have demonstrated that, among rational and selfish agents, a simple competitive reward scheme that reduces the rewards available from already-popular choices can, in theory, return a group to the accuracy implied by the Cordorcet Jury Theorem. This result depends on the assumption that the environmental information received by the agents is truly independent and is not systematically wrong, but effectively balances the expected gains of following social information by choosing the more popular option, and so prevents the information cascades that limit the collective accuracy of sequential decision-making. Under such a reward scheme, and within the assumptions of the model used here, agent's decisions become independent, and depend only on their private information. This increases the probability that the majority of the group will make the correct choice, albeit at the cost of making each individual somewhat less accurate on average. This paper has derived the optimal form and magnitude of this competition in the context of a model in which an agent observes the full sequence of previous decisions, but because agents' decisions become independent under the optimal competition it would retain the same form if agents instead observed simplified aggregate statistics regarding how many agents had made each choice \cite{mann2021optimal}. As such it is applicable across a wide range of domains where the nature of social information may vary.

Introducing competition that penalises agents for following popular choices is an established mechanism for motivating agents to make decisions that improve collective accuracy by reducing the correlation between different decision-makers \cite{hong, mann2017optimal}, and is an important feature of markets as a forecasting mechanism, whether explicitly prediction markets \cite{wolfers2006pmi}, betting exchanges or financial markets. In this paper I have shown that competition can also fulfil this role in a sequential decision-making context where agents can observe the choices made by all those who decide before them and utilise that information in their own decision-making. While the optimal level of competition is unlikely to be known a priori for any given decision or decision-making system, sensitivity analysis shows that introducing a small degree of competition typically improves upon performance from binary rewards alone; in an adaptive system competitive pressure can thus be gradually raised to determine optimal performance. Except in very difficult decisions ($q \simeq 0.5$), competitive rewards are relatively forgiving to miscalibration, providing improvements on binary rewards across a wide scale of competition strengths. 

The theoretical efficacy of competitive rewards raises the possibility that such incentives could be used to improve collective accuracy across a range of real world contexts. For example, the collective judgement of the scientific community  (as reflected in majority expressed opinion) on issues where there is significant uncertainty could potentially be improved by systematically assigning greater rewards to those later proved correct when fewer others also expressed that opinion; these rewards might be in the form of promotions, research funding or simply scientific reputation. To some degree such competitive rewards already feature in many communities, and many scientists, economists and political pundits have made their reputation by advocating for a minority viewpoint that was later proved correct: a notable example is the case of Barry Marshall and Robin Warren, who won the 2005 Nobel Prize in Physiology or Medicine for their discovery of the link between \emph{H. pylori} and stomach ulcers \cite{marshall1985atf}. However, other pressures that incentivise social and professional conformity are also common, such as what Irving Janis termed `Groupthink' \cite{janis1983groupthink} -- the tendency to excessively value consensus with other group members. Conformity may also be imposed by systemic factors such as needing to convince others that your ideas are plausible before they can be explored \cite{gross2021wep}. 

From an external point of view, the results of this study suggest we should assign greater credibility to the collective wisdom of communities where such competitive rewards are the norm, motivating both accuracy and independence. Conversely, the collective wisdom of communities characterised by strong social norms of conformity (effectively negative competition) should be assigned lower credibility. Notably, although competitive, some communities such as political punditry rarely demand or reward specific, falsifiable predictions \cite{tetlock2016superforecasting}; for competitive rewards to drive collective accuracy there must be penalties (or lack of reward) for inaccurate predictions, otherwise individuals are simply motivated to identify and state a unique opinion without regard for its accuracy. The optimal reward structure identified here requires that rewards are still contingent on accuracy.

Where the collective accuracy of group decisions (at least as expressed via majority voting) is highly desirable, we should seek to reduce pressures that induce conformity and introduce competitive rewards that motivate more independent judgements. However, as well as potentially causing social friction (if norms of social conformity are violated), this also comes at the cost of a likely reduction in individual accuracy. While the group may be more accurate, more individuals will be wrong. This highlights that systems of collective decision-making that aim for collective accuracy must not only seek to be tolerant of conflicting views, but must also tolerate a greater level of individual decision-making failure.

To be maximally effective such competitive rewards need to be predictable. Agents should be able to either rationally adjust their choices in the light of known reward schemes, or such rewards should be consistent enough to allow adaptation by reinforcement learning -- a process which may take time to affect behavioural change \cite{burton2015payoff, burton2021payoff, burton2021selfish}. Rewards should also be well-tuned to the specific context (particularly the level of individual certainty, but also the size of the community). This suggests that competitive reward structures should be made more explicit, and calibrated through systematic trial and error in a particular community. The model presented here is theoretical and assumes that agents are either well informed about potential rewards and respond rationally, or reliably adapt via reinforcement based on experience. Such assumptions, and the efficacy of innovative collective decision-making systems are are ultimately empirical questions that must be tested experimentally. 

Finally, this paper has considered problems in which agents seek to ascertain the answer to a question of external empirical fact, such as whether it will rain tomorrow, or which of two teams will a sporting event. It should be noted that some attempts to leverage the predictive power of social information instead focus on questions where the answer is endogenous to the community from which that social information is drawn. An example is the use of social media to predict which movies will attract large box office returns \cite{asur2010predicting}, since presumably the commenters on social media represent a sample of potential movie-goers. In these cases there is likely to be more value in simply measuring the aggregate opinion of individuals, since the expression of interest in a movie is itself a predictor of attendance, regardless of whether that interest is itself socially driven. 

\subsection*{Acknowledgments}
This work was supported by UK Research and Innovation Future Leaders Fellowship MR/S032525/1.

\bibliographystyle{ieeetr}

\begin{thebibliography}{10}

\bibitem{de2014essai}
N.~De~Condorcet, {\em Essai sur l'application de l'analyse {\`a} la
  probabilit{\'e} des d{\'e}cisions rendues {\`a} la pluralit{\'e} des voix}.
\newblock Paris, 1785.

\bibitem{galton1907vox}
F.~Galton, ``Vox populi,'' {\em Nature}, vol.~75, no.~7, pp.~450--451, 1907.

\bibitem{surowiecki2005wisdom}
J.~Surowiecki, {\em {The Wisdom of Crowds}}.
\newblock Random House LLC, 2005.

\bibitem{list2001epistemic}
C.~List and R.~E. Goodin, ``{Epistemic democracy: Generalizing the Condorcet
  jury theorem},'' {\em Journal of Political Philosophy}, vol.~9, no.~3,
  p.~277306, 2001.

\bibitem{landemore2012democratic}
H.~Landemore, {\em Democratic reason}.
\newblock Princeton University Press, 2012.

\bibitem{klein2011harvest}
M.~Klein, ``{How to harvest collective wisdom on complex problems: An
  introduction to the MIT deliberatorium},'' {\em Center for Collective
  Intelligence working paper}, 2011.

\bibitem{boland1989majority}
P.~J. Boland, ``{Majority systems and the Condorcet jury theorem},'' {\em
  Journal of the Royal Statistical Society: Series D}, vol.~38, no.~3,
  pp.~181--189, 1989.

\bibitem{list2004democracy}
C.~List, ``Democracy in animal groups: a political science perspective,'' {\em
  Trends in Ecology and Evolution}, vol.~19, no.~4, pp.~168--169, 2004.

\bibitem{king2007use}
A.~J. King and G.~Cowlishaw, ``When to use social information: the advantage of
  large group size in individual decision making,'' {\em Biology Letters},
  vol.~3, no.~2, pp.~137--139, 2007.

\bibitem{bikhchandani1992theory}
S.~Bikhchandani, D.~Hirshleifer, and I.~Welch, ``A theory of fads, fashion,
  custom, and cultural change as informational cascades,'' {\em Journal of
  Political Economy}, vol.~100, no.~5, pp.~992--1026, 1992.

\bibitem{faria2010coll}
J.~Faria, S.~Krause, and J.~Krause, ``Collective behavior in road crossing
  pedestrians: the role of social information,'' {\em Behavioral Ecology},
  vol.~21, no.~6, pp.~1236--1242, 2010.

\bibitem{gallup2012visual}
A.~C. Gallup, J.~J. Hale, D.~J.~T. Sumpter, S.~Garnier, A.~Kacelnik, J.~R.
  Krebs, and I.~D. Couzin, ``Visual attention and the acquisition of
  information in human crowds,'' {\em Proceedings of the National Academy of
  Sciences}, vol.~109, no.~19, pp.~7245--7250, 2012.

\bibitem{mann2013tdo}
R.~P. Mann, J.~Faria, D.~J.~T. Sumpter, and J.~Krause, ``The dynamics of
  audience applause,'' {\em Journal of the Royal Society Interface}, vol.~10, p.~20130466, 2013.

\bibitem{sumpter2009quo}
D.~J.~T. Sumpter and S.~C. Pratt, ``Quorum responses and consensus decision
  making,'' {\em Philosophical Transactions of the Royal Society B: Biological
  Sciences}, vol.~364, no.~1518, pp.~743--753, 2009.

\bibitem{anderson1997information}
L.~R. Anderson and C.~A. Holt, ``Information cascades in the laboratory,'' {\em
  The American Economic Review}, pp.~847--862, 1997.

\bibitem{mann2018cdm}
R.~P. Mann, ``Collective decision making by rational individuals,'' {\em
  Proceedings of the National Academy of Sciences}, vol.~115, no.~44,
  pp.~E10387--E10396, 2018.

\bibitem{tump2020wise}
A.~N. Tump, T.~J. Pleskac, and R.~H. Kurvers, ``Wise or mad crowds? the
  cognitive mechanisms underlying information cascades,'' {\em Science
  Advances}, vol.~6, no.~29, p.~eabb0266, 2020.

\bibitem{mackay2012extraordinary}
C.~Mackay, {\em Extraordinary popular delusions and the madness of crowds}.
\newblock Farrar, Straus and Giroux: New York, 1841.

\bibitem{lorenz2011hsi}
J.~Lorenz, H.~Rauhut, F.~Schweitzer, and D.~Helbing, ``How social influence can
  undermine the wisdom of crowd effect,'' {\em Proceedings of the National
  Academy of Sciences}, vol.~108, no.~22, pp.~9020--9025, 2011.

\bibitem{hong}
L.~Hong, S.~E. Page, and M.~Riolo, ``Incentives, information, and emergent
  collective accuracy,'' {\em Managerial and Decision Economics}, vol.~33,
  no.~5-6, pp.~323--334, 2012.

\bibitem{mann2017optimal}
R.~P. Mann and D.~Helbing, ``Optimal incentives for collective intelligence,''
  {\em Proceedings of the National Academy of Sciences}, vol.~114, no.~20,
  pp.~5077--5082, 2017.

\bibitem{perez2011cab}
A.~P{\'e}rez-Escudero and G.~G. De~Polavieja, ``{Collective Animal Behavior
  from Bayesian Estimation and Probability Matching},'' {\em PLoS Computational
  Biology}, vol.~7, no.~11, p.~e1002282, 2011.

\bibitem{arganda2013acr}
S.~Arganda, A.~P{\'e}rez-Escudero, and G.~G. De~Polavieja, ``A common rule for
  decision-making in animal collectives across species,'' {\em Proceedings of
  the National Academy of Sciences}, vol.~109, pp.~20508--20513, 2012.

\bibitem{mann2020collective}
R.~P. Mann, ``Collective decision-making by rational agents with differing
  preferences,'' {\em Proceedings of the National Academy of Sciences},
  vol.~117, no.~19, pp.~10388--10396, 2020.

\bibitem{aumann1976agreeing}
R.~J. Aumann, ``Agreeing to disagree,'' {\em The Annals of Statistics}, vol.~4,
  no.~6, pp.~1236--1239, 1976.

\bibitem{torney2015social}
C.~J. Torney, T.~Lorenzi, I.~D. Couzin, and S.~A. Levin, ``Social information
  use and the evolution of unresponsiveness in collective systems,'' {\em
  Journal of the Royal Society Interface}, vol.~12, no.~103, p.~20140893, 2015.

\bibitem{mann2021optimal}
R.~P. Mann, ``Optimal use of simplified social information in sequential
  decision-making,'' {\em Journal of the Royal Society Interface}, vol.~18,
  no.~179, p.~20210082, 2021.

\bibitem{wolfers2006pmi}
J.~Wolfers and E.~Zitzewitz, ``Prediction markets in theory and practice,''
  {\em Journal of Economic Perspectives}, vol.~18, pp.~107--126, 2004.

\bibitem{marshall1985atf}
B.~J. Marshall, J.~A. Armstrong, D.~B. McGechie, and R.~J. Clancy, ``{Attempt
  to fulfil Koch's postulates for pyloric Campylobacter},'' {\em Medical
  Journal of Australia}, vol.~142, no.~8, pp.~436--439, 1985.

\bibitem{janis1983groupthink}
I.~L. Janis, {\em Groupthink}.
\newblock Houghton Mifflin Boston, 1983.

\bibitem{gross2021wep}
K.~Gross and C.~T. Bergstrom, ``Why ex post peer review encourages high-risk
  research while ex ante review discourages it,'' {\em Proceedings of the
  National Academy of Sciences}, vol.~118, no.~51, p.~e2111615118, 2021.

\bibitem{tetlock2016superforecasting}
P.~E. Tetlock and D.~Gardner, {\em Superforecasting: The art and science of
  prediction}.
\newblock Random House, 2016.

\bibitem{burton2015payoff}
M.~N. Burton-Chellew, H.~H. Nax, and S.~A. West, ``Payoff-based learning
  explains the decline in cooperation in public goods games,'' {\em Proceedings
  of the Royal Society of London B: Biological Sciences}, vol.~282, no.~1801,
  p.~20142678, 2015.

\bibitem{burton2021payoff}
M.~N. Burton-Chellew and S.~A. West, ``Payoff-based learning best explains the
  rate of decline in cooperation across 237 public-goods games,'' {\em Nature
  Human Behaviour}, pp.~1--9, 2021.

\bibitem{burton2021selfish}
M.~N. Burton-Chellew and C.~Gu{\'e}rin, ``Selfish learning is more important
  than fair-minded conditional cooperation in public-goods games,'' {\em
  Available at SSRN 3982907}, 2021.

\bibitem{asur2010predicting}
S.~Asur and B.~A. Huberman, ``Predicting the future with social media,'' in
  {\em 2010 IEEE/WIC/ACM international conference on web intelligence and
  intelligent agent technology}, vol.~1, pp.~492--499, IEEE, 2010.

\end{thebibliography}

\end{document}